%% file: lowcontrast.tex
\documentclass{article}
\usepackage{amsmath, amsthm, amssymb, amsfonts}
\usepackage{graphicx}
\usepackage{subfigure}
\usepackage{pgfplots}
\pgfplotsset{compat=1.4}
  \usepackage{paralist}
 \usepackage[colorlinks=true]{hyperref}
\hypersetup{urlcolor=blue, citecolor=red}

\newtheorem{theorem}{Theorem}[section]

\newtheorem{proposition}{Proposition}

\newtheorem{definition}[theorem]{Definition}
\newtheorem{remark}{Remark}
\newcommand{\ep}{\varepsilon}

\title{Motion of discrete interfaces\\ in low-contrast periodic media}

\author{}
\date{}

\begin{document}
\maketitle

\centerline{\scshape Giovanni Scilla }
\medskip
{\footnotesize
 \centerline{Dipartimento di Matematica `G. Castelnuovo'}
   \centerline{`Sapienza' Universit\`a di Roma}
   \centerline{piazzale Aldo Moro 5, 00185 Roma (Italy)}
} 
%
%


\begin{abstract}
We study the motion of discrete interfaces driven by ferromagnetic interactions in a two-dimensional low-contrast periodic environment, by coupling the minimizing movements approach by Almgren, Taylor and Wang and a discrete-to-continuum analysis. As in a recent paper by Braides and Scilla dealing with high-contrast periodic media, we give an example showing that in general the effective motion does not depend only on the $\Gamma$\hbox{-}limit, but also on geometrical features that are not detected in the static description. We show that there exists a critical value $\widetilde{\delta}$ of the contrast parameter $\delta$ above which the discrete motion is constrained and coincides with the high-contrast case. If $\delta<\widetilde{\delta}$ we have a new pinning threshold and a new effective velocity both depending on $\delta$. We also consider the case of non-uniform inclusions distributed into periodic uniform layers.
\end{abstract}

\section{Introduction}

In this paper we study a problem of homogenization for a discrete crystalline flow. The analysis will be carried over by using the minimizing-movement scheme of Almgren, Taylor and Wang \cite{ATW83}. This consists in introducing a time scale $\tau$, iteratively defining a sequence of sets $E^\tau_k$ as minimizers of

\begin{equation}
\min \Bigl\{ P(E)+{\frac{1}{\tau}} D(E, E^\tau_{k-1})\Bigr\} ,
\end{equation}
\\
where $P$ is a perimeter energy and $D$ is a distance-type energy between sets, and $E^\tau_0$
is a given initial datum, and subsequently computing a time-continuous limit $E(t)$ of $\{E^\tau_k\}$ as $\tau\to0$, which defines the desired geometric motion related to the energy~$P$.

In recent papers by Braides, Gelli and Novaga \cite{BGN} and Braides and Scilla \cite{BraSci13}, the Almgren-Taylor-Wang approach has been used coupled to a homogenization procedure. In this case the perimeters and the distances depend on a small parameter $\ep$ (interpreted as a space scale), and consequently, after introducing a time scale $\tau$, the time-discrete motions are the $E^{\tau,\ep}_k$ defined iteratively by

\begin{equation}
E^{\tau,\ep}_k \hbox{ is a minimizer of } \min \Bigl\{ P_\ep(E)+{\frac{1}{\tau}} D_\ep(E, E^{\tau,\ep}_{k-1})\Bigr\}.
\end{equation}

The time-continuous limit $E(t)$ of $\{E^{\tau,\ep}_k\}$ then may depend how mutually $\ep$ and $\tau$ tend to $0$ (see Braides~\cite{Bra13}). In particular, if we have a large number of local minimizers then the limit motion will be pinned if $\tau<\!\!<\ep$ suitably fast (in a sense, we can pass to the limit in $\tau$ first, and then apply the Almgren-Taylor-Wang approach, which clearly gives pinning when the initial data are local minimizers). On the contrary, if $\ep<\!\!<\tau$ fast enough
and $P_\ep$ $\Gamma$-converge to a limit perimeter $P$ (which is always the case by compactness) then the limit $E$ will be
the evolution related to the limit $P$ (again, in a sense, in this case we can pass to the limit in $\ep$ first).

\begin{figure}[htbp]
\centering
\def\svgwidth{100pt}
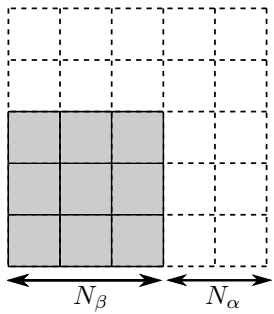
\caption{The periodicity cell.}\label{fig:0}
\end{figure}

In \cite{BraSci13} the energies $P_\varepsilon$ are \emph{inhomogeneous ferromagnetic energies} defined on subsets $E\subset \varepsilon\mathbb{Z}^2$, of the form

\begin{equation*}
P_\ep(E)= {\frac{1}{2}}\ep\,\sum\{c_{ij}: i, j\in \mathbb{Z}^2,\ep i\in E, \ep j\not\in E, \ |i-j|=1\},
\end{equation*}
\\
(we use the notation $\sum\{x_a:a\in A\}=\sum_{a\in A}x_a$) where the coefficients $c_{ij}$ equal $\alpha>0$ except for some well-separated periodic square
inclusions of size $N_\beta$ where $c_{ij}=\beta>\alpha$. The periodicity cell is pictured in Fig.~\ref{fig:0}, where continuous lines represent $\beta$-bonds, dashed lines $\alpha$-bonds.
These inclusions are not energetically favorable (\emph{high\hbox{-}contrast} medium) and they can be neglected in the computation of the $\Gamma$-limit, which is the \emph{crystalline perimeter}

$$
P(E)=\alpha\int_{\partial E}\|\nu\|_1d\mathcal{H}^1,
$$
\\
where $\nu$ is the normal to $\partial E$ and $\|(\nu_1,\nu_2)\|_1=|\nu_1|+|\nu_2|$.

The {flat flow} of this perimeter is the motion by crystalline curvature in dimension two described by Taylor \cite{Ta}. In the case of initial datum a coordinate rectangle, the evolution by crystalline curvature is a rectangle with the same centre and sides of lengths $L_1, L_2$ governed by the system of ordinary differential equations

$$
\begin{cases}\displaystyle \dot L_1= -{\frac{4\alpha}{L_2}}\cr\cr
\displaystyle \dot L_2= -{\frac{4\alpha}{L_1}}.\end{cases}
$$

In \cite{BraSci13} all possible evolutions have been characterized as $\ep, \tau\to 0$ showing that the relevant mutual scale is when $\tau/\ep\to\gamma$. In the case of initial datum a coordinate rectangle the resulting evolution is still a rectangle. In the case of a unique evolution the side lengths $L_1(t), L_2(t)$ of this rectangle are governed by a system of `degenerate' ordinary differential equations

\begin{equation}
\begin{cases}\displaystyle \dot L_1= -{\frac{2}{\gamma}}f\left({\frac{\gamma}{L_2}}\right)\\ \cr
\displaystyle \dot L_2= -{\frac{2}{\gamma}}f\left({\frac{\gamma}{L_1}}\right),
\end{cases}
\label{systemdeg}
\end{equation}
\\
where the \emph{effective velocity function} $f$, obtained as solution of a one-dimensional homogenization problem, is  a locally constant function on compact subsets of $(0,+\infty)$  which depends on $\alpha$, the period and size of the inclusions but not on $\gamma$ (neither on the value $\beta$). This function has been computed, by means of algebraic formulas, in the simpler cases $N_\beta=1$ and $N_\beta=2$, with varying $N_\alpha$. In particular, if $N_\alpha=N_\beta=1$, then the velocity function is given by

\begin{equation*}
\overline{f}(\gamma/L)=2\left\lfloor\frac{\alpha \gamma}{L}+\frac{1}{4}\right\rfloor,
\end{equation*}
\\
while in case of no inclusion (i.e., $\alpha=\beta$), it is given by

\begin{equation*}
\widetilde{f}(\gamma/L)=\left\lfloor\frac{2\alpha \gamma}{L}\right\rfloor.
\end{equation*}
\\
The dependence on the inclusions gives the \emph{pinning threshold} (i.e., the critical value of the side length above which it is energetically not favorable for a side to move)

$$
\overline L={\frac{4\alpha\gamma}{2+N_\beta}}
$$
\\
depending on the size  of the inclusion $N_\beta$.

The inclusions can be considered as ``obstacles'' that can be bypassed when computing minimizers of $P_\ep$; however their presence is felt in the minimizing-movement procedure since they may influence the choice of $E^{\tau,\ep}_k$ through the interplay between the distance and perimeter terms. As a result the motion can be either decelerated or accelerated with respect to the homogeneous case (i.e., the case $\alpha=\beta$ described in \cite{BGN}).\\

Scope of this work is to give another example showing that the periodic microstructure can affect the limit evolution without changing the $\Gamma$-limit. To this end we perform a multi-scale analysis by introducing a \emph{contrast parameter} $\delta_\ep$ and considering a \emph{low-contrast} medium, that is a periodic mixture of two homogeneous materials whose propagating properties are close to each other (see e.g. \cite{Milton}). One of them can be considered as a fixed background medium (described by $\alpha$-connections) and the other as a small (vanishing) perturbation from that one, that is with $\beta=\beta_\ep=\beta(\ep)$ and $\beta_\ep-\alpha=\delta_\ep\to0$ as $\ep\to0$. With the same notation as in \cite{BraSci13} we restrict ourselves to the case $N_\alpha=N_\beta=1$ (see Fig.~\ref{fig:0}); despite of its simplicity, the choice of this particular geometry will suffice to show new features of the motion. The main result is the existence of a threshold value of the contrast parameter below which we have a new homogenized effective velocity, which takes into account the propagation velocities in both the connections $\alpha$ and $\beta$; above this threshold, instead, it is independent of the value of $\beta$ and the motion is obtained by considering only the $\alpha$\hbox{-}connections. The dependence of the effective properties on microstructure in low-contrast periodic media has been recently investigated for various physical problems (see e.g. \cite{Con}).

\begin{figure}
\centering
\def\svgwidth{250pt}
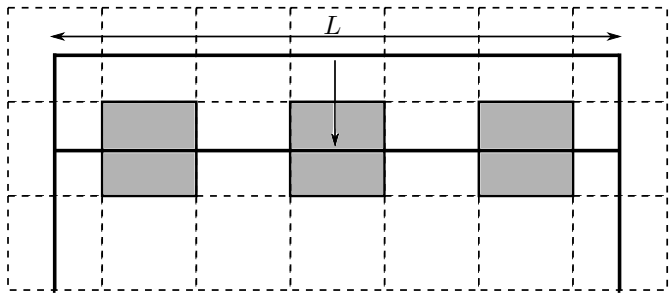
\caption{Displacement of a side from $\alpha$-connections to $\beta$-connections.}
\label{heu}
\end{figure}

We first determine the correct scaling for $\delta_\ep$ to have that a straight interface may stay on $\beta$-connections. To this end we consider a coordinate rectangle whose sides intersect only $\alpha$-bonds (\emph{$\alpha$-type rectangle}), we write the variation of the energy $\mathcal{F}_{\ep,\gamma\ep}^{\alpha,\beta_\ep}$ (\ref{newenergy}) associated to the displacement by $\ep$ of the upper horizontal side of length $L$ (see Fig.~\ref{heu}) and we impose it to be zero. We have that

\begin{equation*}
-2\alpha\ep+(\beta_\ep-\alpha)L+\frac{cL}{\gamma}\ep=-2\alpha\ep+\delta_\ep L+\frac{cL}{\gamma}\ep = 0,
\end{equation*}
\\
where $c=c(L)$ is a constant depending on $L$, from which we deduce that

\begin{equation*}
\delta_\ep=\left(\frac{2\alpha}{L}-\frac{c}{\gamma}\right)\ep=O(\ep) \quad \text{as }\ep\to0.
\end{equation*}
\\
This heuristic computation suggests that the correct scaling is
\begin{equation*}
\beta_\ep-\alpha=\delta_\ep=\delta\ep
\end{equation*}
for some constant $\delta>0$.

As in \cite{BraSci13} we will assume that $\tau=\gamma\ep$ and we will restrict the description of the motion to the case of initial data coordinate rectangles, since all other cases can be reduced to the study of this one. The evolution of a coordinate rectangle by minimization of the energy is again a coordinate rectangle.
We will show that there exists a threshold $\widetilde{\delta}=\frac{1}{2\gamma}$ such that if $\delta<\widetilde{\delta}$ (subcritical regime) then the evolution is a rectangle that may have some \emph{$\beta$-type side} (that is, a side intersecting only $\beta$\hbox{-}connections), while if $\delta\geq\widetilde{\delta}$ (supercritical regime) the $\beta$\hbox{-}connections are avoided as in the case $\beta>\alpha$. Note that this result gives information also for more general choices of the vanishing rate of $\delta_\ep$: if $\delta_\ep<\!<\!\ep$ we reduce to the subcritical case, while if $\delta_\ep >\!>\!\ep$ we are in the supercritical case. The limit motion can still be described through a system of degenerate ordinary differential equations as in (\ref{systemdeg}) with a new effective velocity function $f=f_\delta$ depending on $\delta$. We also have a new effective pinning threshold given by

\begin{equation*}
\overline{L}_\delta=\max\left\{\frac{2\alpha\gamma}{\delta\gamma+1},\frac{4}{3}\alpha\gamma\right\}.
\end{equation*}

The paper is organized as follows. In Section~\ref{setting} we define all the energies that we will consider. We then formulate the discrete-in-time scheme analogous to the Almgren, Taylor and Wang approach. Section~\ref{rectangle} contains the description of the convergence of the discrete scheme in the case of a rectangular initial set. We show that the minimizers of this scheme are actually rectangles. Subsection~\ref{effveloc} deals with the definition of the effective velocity of a side by means of a homogenization formula resulting from a one-dimensional `oscillation-optimization' problem. This velocity can be expressed uniquely (up possibly to a discrete set of values) as a function of the ratio of $\gamma$ and the side length, and of $\delta$ (Definition~\ref{effvel}). Subsection~\ref{newpinning} contains the computation of the effective pinning threshold, showing that it is affected by microstructure because it also depends on the parameter $\delta$. In Subsection~\ref{computation} we compute explicitly the velocity function showing a comparison with the homogeneous case $\alpha=\beta$ and the high-contrast case $\beta>\alpha$. The description of the homogenized limit motion is contained in Subsection~\ref{limitmotion}.  Section~\ref{periodic2} deals with the case of non-uniform inclusions distributed into periodic uniform layers.

\section{Notation and setting of the problem}\label{setting}

If $x=(x_1,x_2)\in\mathbb{R}^2$ we set $\|x\|_1=|x_1|+|x_2|$ and $\|x\|_\infty=\max\{|x_1|,|x_2|\}$. If $A$ is a Lebesgue\hbox{-}measurable set we denote by $|A|$ its two-dimensional Lebesgue measure. The symmetric difference of $A$ and $B$ is denoted by $A\triangle B$, their Hausdorff distance by $\text{d}_\mathcal{H}(A,B)$. If $E$ is a set of finite perimeter then $\partial^*E$ is its reduced boundary (see, for example \cite{Bra98}). The measure-theoretical inner normal to $E$ at a point $x$ in $\partial^*E$ is denoted by $\nu=\nu_E(x)$.

\subsection{Inhomogeneous `low-contrast' ferromagnetic energies}\label{inhom}

\begin{figure}
\center
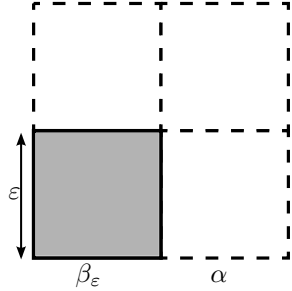
\caption{The periodicity cell.}
\label{fig1}
\end{figure}

The energies we consider are interfacial energies defined in an inhomogeneous low-contrast environment as follows. Let $\alpha,\delta>0$ and we fix $\ep>0$. We consider $2\ep$-periodic coefficients $c_{ij}^\ep$ indexed on \emph{nearest-neighbors} of $\ep\mathbb{Z}^2$ (i.e., $i,j\in\ep\mathbb{Z}^2$ with $|i-j|=\ep$) defined for $i,j$ such that

\begin{equation*}
0\leq \frac{i_1+j_1}{2},\frac{i_2+j_2}{2}<2\ep
\end{equation*}
\\
by
\begin{equation}
c_{ij}^\ep=
\begin{cases}
\beta_\ep=\alpha+\delta\ep, &\quad \text{if }0\leq \displaystyle\frac{i_1+j_1}{2},\frac{i_2+j_2}{2}\leq\ep\\
\alpha   & \quad\text{otherwise.}
\end{cases}
\label{coeff}
\end{equation}
\\
These coefficients label the bonds between points in $\ep\mathbb{Z}^2$, so that they describe a matrix of $\alpha$-bonds with $2\ep$-periodic inclusions of $\beta$-bonds grouped in squares of side length $\ep$. The periodicity cell is pictured in Fig.~\ref{fig1}. Here the continuous lines represent $\beta$-bonds while the dashed lines the $\alpha$ ones.

Correspondingly, to coefficients (\ref{coeff}) we associate the energy defined on subsets $\mathcal{I}$ of $\ep\mathbb{Z}^2$ by
\begin{equation}
\text{P}^{\alpha,\beta_\ep}_\ep(\mathcal{I})=\sum_{i,j}\left\{\ep c_{ij}^\ep: |i-j|=\ep, i\in\mathcal{I},j\in\ep\mathbb{Z}^2\backslash\mathcal{I}\right\},
\label{energy}
\end{equation}
where, as remarked in the Introduction, we use the notation $\sum\{x_a:a\in A\}=\sum_{a\in A}x_a$.

In order to study the continuous limit as $\ep\to0$ of these energies, we will identify each subset of $\ep\mathbb{Z}^2$ with a measurable subset of $\mathbb{R}^2$ in such a way that equi-boundedness of the energies implies pre-compactness of such sets in the sense of the sets of finite perimeter. This identification is as follows: if $\ep>0$ and $i\in\ep\mathbb{Z}^2$, we denote by $Q_\ep(i)$ the closed coordinate square with side length $\ep$ and centered in $i$. To a set of indices $\mathcal{I}\subset\ep\mathbb{Z}^2$ we associate the set
\begin{equation}
E_\mathcal{I}=\bigcup_{i\in\mathcal{I}}Q_\ep(i).
\end{equation}
The space of \emph{admissible sets} related to indices in the two-dimensional square lattice is then defined by
\begin{equation*}
\mathcal{D}_\ep:=\{E\subseteq\mathbb{R}^2:\quad E=E_\mathcal{I}\text{ for some }\mathcal{I}\subseteq\ep\mathbb{Z}^2\}.
\end{equation*}
For each $E=E_\mathcal{I}\in\mathcal{D}_\ep$ we denote
\begin{equation*}
\text{P}_\ep^{\alpha,\beta_\ep}(E)=\text{P}_\ep^{\alpha,\beta_\ep}(\mathcal{I}).
\end{equation*}
We note that
\begin{equation}
\text{P}_\ep^{\alpha,\beta_\ep}(E)\geq\ep\alpha\left\{(i,j): |i-j|=\ep,  i\in\mathcal{I},j\in\ep\mathbb{Z}^2\backslash\mathcal{I}\right\}=\alpha\mathcal{H}^1(\partial E),
\end{equation}
which shows that sequences of sets $E_\ep$ with $\sup_\ep\text{P}^{\alpha,\beta_\ep}_\ep(E_\ep)<+\infty$ are pre\hbox{-}compact with respect to the local $L^1$-convergence in $\mathbb{R}^2$ of their characteristic functions and their limits are sets of finite perimeter in $\mathbb{R}^2$.
Hence, this defines a meaningful convergence with respect to which compute the $\Gamma$\hbox{-}limit of $\text{P}^{\alpha,\beta_\ep}_\ep$ as $\ep\to0$.
The energies $\text{P}_\ep^{\alpha,\beta_\ep}$ defined by (\ref{energy}) $\Gamma$\hbox{-}converge, as $\ep\to0$, to the anisotropic cristalline perimeter functional
\begin{equation*}
\text{P}^\alpha(E)=\alpha\int_{\partial^*E}\|\nu\|_1\,d\mathcal{H}^1.
\end{equation*}
This can be shown with an analogous computation as in Braides-Scilla \cite{BraSci13}.

\subsection{A discrete distance}

For $\mathcal{I}\subset\ep\mathbb{Z}^2$ we define the \emph{discrete $\ell^{\infty}$-distance} from $\partial\mathcal{I}$ as
\begin{equation*}
d_\infty^\ep(i,\partial\mathcal{I})=
\begin{cases}
\inf\{\|i-j\|_\infty:j\in\mathcal{I}\}&\text{if $i\not\in\mathcal{I}$}\\
\inf\{\|i-j\|_\infty:j\in\ep\mathbb{Z}^2\backslash\mathcal{I}\}&\text{if $i\in\mathcal{I}$}.
\end{cases}
\end{equation*}
Note that we have $d_\infty^\ep(i,\partial\mathcal{I})=d_\infty(i,\partial E_\mathcal{I})+\displaystyle\frac{\ep}{2}$, where $d_\infty$ denotes the usual $\ell^\infty$-distance. The distance can be extended to all $\mathbb{R}^2\backslash\partial E_\mathcal{I}$ by setting
\begin{equation*}
d_\infty^\ep(x,\partial\mathcal{I})=d_\infty^\ep(i,\partial\mathcal{I})\quad \text{if }x\in Q_\ep(i).
\end{equation*}
In the following we will directly work with $E\in\mathcal{D}_\ep$, so that the distance can be equivalently defined by
\begin{equation*}
d_\infty^\ep(x,\partial E)=d_\infty(i,\partial E)+\frac{\ep}{2},\quad \text{if }x\in Q_\ep(i).
\end{equation*}
Note that this is well defined as a measurable function, since its definition is unique outside the union of the boundaries of the squares $Q_\ep$ (that are a negligible set).

\subsection{Minimization scheme}\label{timemin}
We fix a time step $\tau>0$ and introduce a discrete motion with underlying time step $\tau$ obtained by successive minimization. At each time step we will minimize an energy $\mathcal{F}_{\ep,\tau}^{\alpha,\beta_\ep}:\mathcal{D}_\ep\times\mathcal{D}_\ep\to\mathbb{R}$ defined as
\begin{equation}
\mathcal{F}_{\ep,\tau}^{\alpha,\beta_\ep}(F,E)= \text{P}_\ep^{\alpha,\beta_\ep}(F)+\frac{1}{\tau}\int_{F\triangle E}d_\infty^\ep(x,\partial E)\,dx.
\label{newenergy}
\end{equation}
Note that the integral can be indeed rewritten as a sum on the set of indices $\ep\mathbb{Z}^2\cap(F\triangle E)$. More precisely, if $\mathcal{I}=E\cap\ep\mathbb{Z}^2, \mathcal{J}=F\cap\ep\mathbb{Z}^2$, then
\begin{equation*}
\begin{split}
\mathcal{F}_{\ep,\tau}^{\alpha,\beta_\ep}(\mathcal{J}, \mathcal{I})&=\text{P}_\ep^{\alpha,\beta_\ep}(\mathcal{J})+\frac{1}{\tau}\sum_{i\in\mathcal{J}\triangle\mathcal{I}}\ep^2d_\infty^\ep(i,\partial \mathcal{I})\\
&=\text{P}_\ep^{\alpha,\beta_\ep}(\mathcal{J})+\frac{1}{\tau}\left(\sum_{i\in\mathcal{I}\backslash\mathcal{J}}\ep^2d_\infty(i,\mathcal{I})+\sum_{i\in\mathcal{I}\backslash\mathcal{J}}\ep^2d_\infty(i,\ep\mathbb{Z}^2\backslash\mathcal{I})\right).
\end{split}
\end{equation*}

Given an initial set $E^0_{\ep}\in\mathcal{D}_\ep$, we define recursively a sequence $E_{\ep,\tau}^k$ in $\mathcal{D}_\ep$ by requiring the following:
\begin{description}
\item[(i)] $E^0_{\ep,\tau}=E^0_{\ep}$;
\item[(ii)] $E_{\ep,\tau}^{k+1}$ is a minimizer of the functional $\mathcal{F}_{\ep,\tau}^{\alpha,\beta_\ep}(\cdot,E_{\ep,\tau}^k)$.
\end{description}
The \emph{discrete flat flow} associated to functionals $\mathcal{F}_{\ep,\tau}^{\alpha,\beta_\ep}$ is thus defined by
\begin{equation}
E_{\ep,\tau}(t)=E_{\ep,\tau}^{\lfloor t\slash\tau\rfloor},\quad t\geq0.
\label{disefo}
\end{equation}
Assuming that the initial data $E^0_{\ep}$ tend, in the Hausdorff sense, to a sufficiently regular set $E_0$, we are interested in identifying the motion described by any converging subsequence of $E_{\ep,\tau}(t)$ as $\ep,\tau\to0$.

As remarked in the Introduction, the interaction between the two discretization parameters, in time and space, plays a relevant role in such a limiting process. More precisely, the limit motion depends strongly on their relative decrease rate to 0. If $\ep\!<\!<\tau$ then we may first let $\ep\to0$, so that $\text{P}_\ep^{\alpha,\beta_\ep}(F)$ can be directly replaced by the limit anisotropic perimeter $\text{P}^\alpha(F)$ and $\frac{1}{\tau}\int_{F\bigtriangleup E}d_\infty^\ep(x,\partial E)\,dx$ by $\frac{1}{\tau}\int_{F\bigtriangleup E}d_\infty(x,\partial E)\,dx$. As a consequence, the approximated flat motions tend to the solution of the continuous ones studied by Almgren and Taylor \cite{AT95}. On the other hand, if $\ep\!>\!>\tau$ then there is no motion and $E_{\ep,\tau}^k\equiv E^0_{\ep}$. Indeed, for any $F\neq E^0_{\ep}$ and for $\tau$ small enough we have
\begin{equation*}
\frac{1}{\tau}\int_{F\bigtriangleup E^0_{\ep}}d_\infty^\ep(x,\partial E^0_{\ep})\,dx\geq c\frac{\ep}{\tau}>\text{P}_\ep^{\alpha,\beta_\ep}(E^0_{\ep}).
\end{equation*}
In this case the limit motion is the constant state $E_0$. The meaningful regime is the intermediate case $\tau\sim\ep$.

\section{Motion of a rectangle}\label{rectangle}

As shown in \cite{BGN} the relevant case is when $\ep$ and $\tau$ are of the same order and the initial data are coordinate rectangles $E_\ep^0$, which will be the content of this section.

We assume that
\begin{equation*}
\tau=\gamma\ep\quad\text{for some }\gamma\in(0,+\infty),
\end{equation*}
and, correspondingly, we omit the dependence on $\tau$ in the notation of
\begin{equation*}
E_\ep^k=E^k_{\ep,\tau}(=E^k_{\ep,\gamma\ep}).
\end{equation*}

\begin{definition}
A side intersecting only $\alpha$-bonds (resp., $\beta$-bonds) will be called an \emph{$\alpha$-type side} (resp., \emph{$\beta$-type side}). A coordinate rectangle whose sides are $\alpha$-type sides will be called an \emph{$\alpha$-type rectangle}.
\end{definition}

The first result is that coordinate rectangles evolve into coordinate rectangles. This result will be more precise in the following. In fact, we will show that if $\delta<\frac{1}{2\gamma}$ then the evolution is a rectangle that may have some $\beta$-type side, while if $\delta\geq\frac{1}{2\gamma}$ it has only $\alpha$-type sides (Proposition~\ref{thresh}).

\begin{proposition}
Let $E_\ep^0\in\mathcal{D}_\ep$ be a coordinate rectangle. For all $\eta>0$, if $F$ is a minimizer for the minimum problem for $\mathcal{F}_{\ep,\tau}^{\alpha,\beta_\ep}(\cdot,E_\ep^k),k\geq0$, the sides of $E_\ep^k$ are larger than $\eta$ and $\ep$ is small enough, then $F$ is a coordinate rectangle.
\end{proposition}
\proof
It will suffice to show it for $F=E_\ep^1$. We subdivide the proof into steps.\\
{\bf Step 1: connectedness of $F$ and $\alpha$-rectangularization.} The connectedness of $F$ can be showed as in Braides, Gelli and Novaga \cite{BGN}, because the microstructure does not affect the argument therein. Now consider the maximal $\alpha$-type rectangle $R^\alpha$ with each side intersecting $F$. As in \cite{BraSci13} we call the set $F'=F\cup R^\alpha$ the \emph{$\alpha$-rectangularization} of $F$. This set is either an $\alpha$-type rectangle (and in this case we conclude) or it has some protrusions intersecting $\beta$-bonds. In both cases $\text{P}^{\alpha,\beta_\ep}_\ep(F')\leq\text{P}^{\alpha,\beta_\ep}_\ep(F)$, and the symmetric difference with $E_\ep^0$ decreases. To justify this, note that the $\alpha$-rectangularization reduces (or leaves unchanged) $\text{P}^{\alpha,\alpha}_\ep$ and it reduces the symmetric difference. Moreover, from this fact we deduce the \emph{a priori} estimate $\text{d}_\mathcal{H}(\partial E_\ep^1,\partial E_\ep^0)\leq c(L)\ep$, where $c(L)$ is a constant depending on the length $L$ of the smaller side of $E_\ep^0$.\\
{\bf Step 2: optimal profiles of protrusions on $\beta$-squares.} Now we describe the form of the optimal profiles of the boundary of $F$ intersecting $\beta$-squares. As noted in \cite{BraSci13}, $F$ contains an $\alpha$-type rectangle  $R^\alpha=[\ep m_1,\ep M_1]\times [\ep m_2,\ep M_2]$ and is contained in the $\alpha$-type rectangle
\begin{equation*}
\widetilde{R}^\alpha=[\ep(m_1-1),\ep(M_1+1)]\times[\ep(m_2-1),\ep(M_2+1)],
\end{equation*}
whose sides exceed the ones of $R^\alpha$ by at most $2\ep$. We will analyze separately the possible profiles of $F$ close to each side of $R^\alpha$; e.g., in the rectangle $[\ep(m_1-1),\ep(M_1+1)]\times[\ep M_2,\ep(M_2+1)]$ (i.e., close to the upper horizontal side of $R^\alpha$).

\begin{figure}
\centering
\def\svgwidth{250pt}
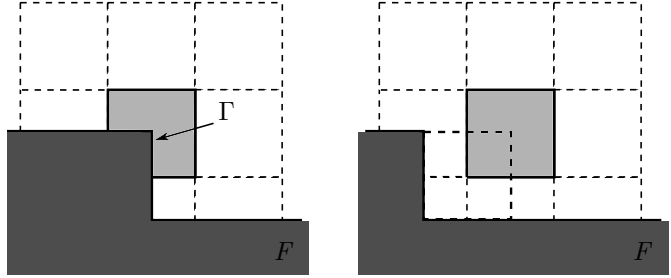
\caption{Removal of an $\ep$-square for $\delta$ `large'.}
\label{rem}
\end{figure}

We first consider the possible behavior of the boundary of $F$ at a single $\beta$\hbox{-}square $Q$, assuming that $Q$ is not one of the two extremal squares. We claim that either $F\cap Q=\emptyset$ or $\partial F\cap Q$ is a horizontal segment. In fact, if a portion $\Gamma$ of $\partial F$ intersects two adjacent sides of $Q$ as in Fig.~\ref{rem}, then we may remove the $\ep$-square whose center is in $Q\cap F$.
In this case, the variation of energy is

\begin{equation}
-2(\beta_\ep-\alpha)\ep+\frac{1}{\gamma}c(L)\ep^2=\left(-2\delta+\frac{1}{\gamma}c(L)\right)\ep^2,
\label{var1}
\end{equation}
\\
which is negative, for $\ep$ small, if and only if $\delta>c(L)/2\gamma$.

\begin{figure}
\centering
\def\svgwidth{250pt}
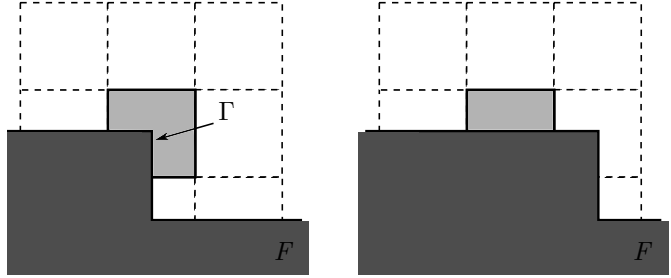
\caption{Adding of an $\ep$-square for $\delta$ `small'.}
\label{add}
\end{figure}

If we add an $\ep$-square as in Fig.~\ref{add}, instead, the variation of the energy is simply

\begin{equation}
-\frac{1}{\gamma}c(L)\ep^2,
\label{var2}
\end{equation}
\\
which is negative. We note that the variation in (\ref{var1}) is strictly less than the one in (\ref{var2}) if and only if $\delta>c(L)/\gamma$.

The same analysis applies to the extremal squares, for which we deduce that $F\cap Q$, if non-empty, is a rectangle with one vertex coinciding with a vertex of $\widetilde{R}^\alpha$.

We now consider the interaction of consecutive $\beta$-squares. Let $Q_1,\dots,Q_K$ a maximal array of consecutive $\beta$-squares with $F\cap Q_k\neq\emptyset$ for $k=1,\dots,K$ and such that $Q_1$ is not a corner square.
\begin{figure}[h]
\centering
\def\svgwidth{350pt}
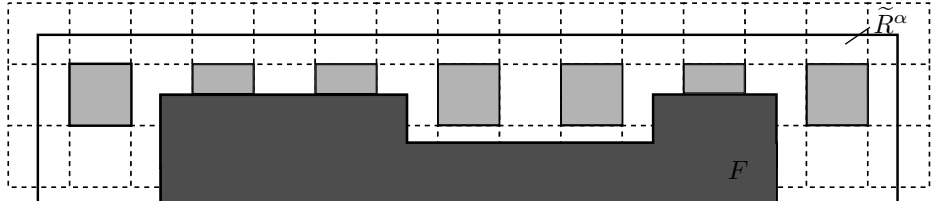
\caption{Interaction of consecutive $\beta$-squares.}
\label{maxim}
\end{figure}
If the subsequent $\beta$-squares $Q_{K+1},\dots,Q_{K+K'}$ are a maximal array not intersecting $F$, and $Q_{K+K'+1},\dots, Q_{K+K'+K''}$ are another maximal array with $F\cap Q_k\neq\emptyset$ for $k=K+K',\dots,K+K'+1$ and such that $Q_{K+K'+K''}$ is not a corner square (see Fig.~\ref{maxim}), then we may replace $F$ by $F\cup R$ (see Fig.~\ref{rep}), where $R$ is the rectangle given by the union of the $\ep$-squares centered at the vertices of the $\beta$-squares $Q_{K+1},\dots,Q_{K+K'}$.
\begin{figure}
\centering
\def\svgwidth{350pt}
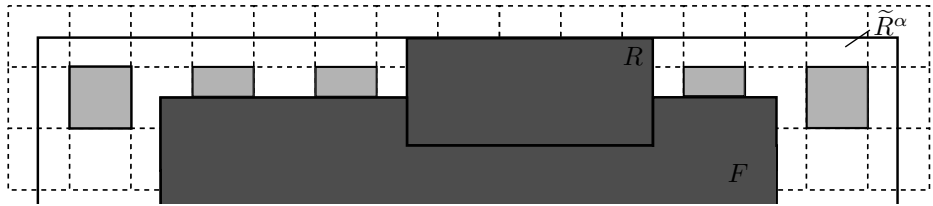
\caption{The new profile after replacing $F$ by $F\cup R$.}
\label{rep}
\end{figure}
This operation leaves unchanged the $\text{P}_\ep^{\alpha,\beta_\ep}$ and reduces the symmetric difference with $E_\ep^0$. We can repeat this procedure for any tern of such arrays. At this point, if we replace $F$ by $F\cup[\ep m_1,\ep M_1]\times[\ep M_2, \ep (M_2+1)]$, this strictly reduces $\text{P}_\ep^{\alpha,\beta_\ep}$ and the symmetric difference (see Fig.~\ref{fin}).

\begin{figure}
\centering
\def\svgwidth{350pt}
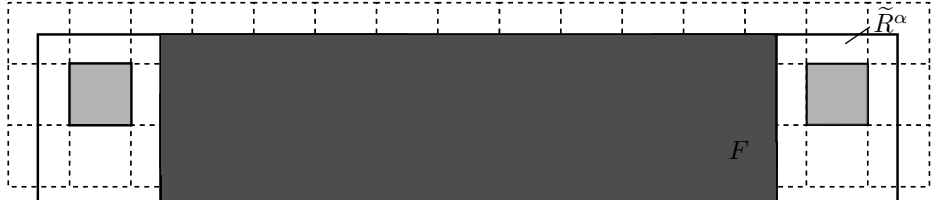
\caption{The new profile after replacing $F$ by $F\cup[\ep m_1,\ep M_1]\times[\ep M_2, \ep (M_2+1)]$.}
\label{fin}
\end{figure}

Note that, if the intersection of $F$ and the left (resp., right) corner square is not empty, then we can consider as a competitor $F\cup[\ep (m_1-1),\ep M_1]\times[\ep M_2, \ep (M_2+1)]$ (resp., $F\cup[\ep m_1,\ep (M_1+1)]\times[\ep M_2, \ep (M_2+1)]$); if $F$ has non empty intersection with both the corner squares, then we consider $F\cup[\ep(m_1-1),\ep(M_1+1)]\times[\ep M_2,\ep(M_2+1)]$.
\begin{figure}[h]
\centering
\def\svgwidth{350pt}
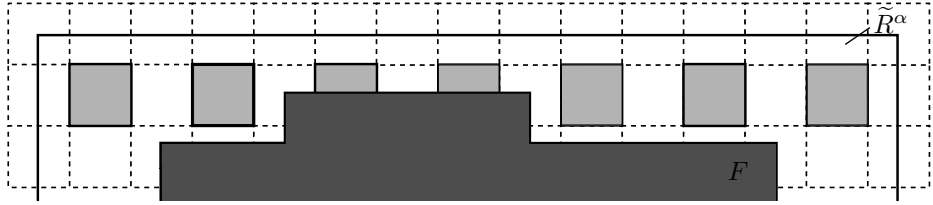
\caption{The case of a single maximal array of intersecting $\beta$-squares.}
\label{only}
\end{figure}
If there exists only one maximal array $Q_1,\dots,Q_K$ and the intersection of $F$ and both the corner squares is empty (see Fig.~\ref{only}), then we may remove all the $\ep$-squares centered at vertices of $Q_1,\dots,Q_K$ and the variation of energy is

\begin{equation}
-2\alpha\ep+2K(\beta_\ep-\alpha)\ep+\frac{1}{\gamma}2c(L)K\ep^2=-2\alpha\ep-2K\delta\ep^2+\frac{1}{\gamma}2c(L)K\ep^2,
\label{comput}
\end{equation}
\\
which is negative for $\ep\leq\frac{\alpha\gamma}{K(c(L)-\delta\gamma)}$ if $\delta<c(L)/\gamma$, any $\ep$ if $\delta\geq c(L)/\gamma$.

\begin{figure}[h]
\centering
\def\svgwidth{350pt}
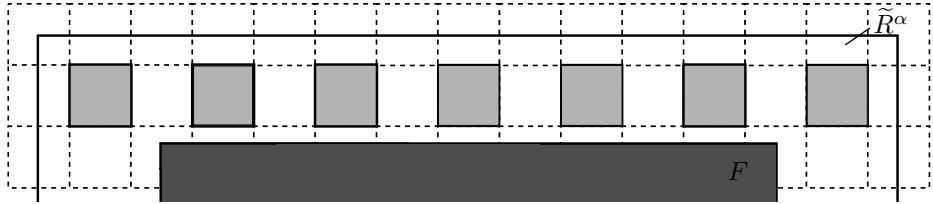
\caption{The profile after removing all the $\ep$-squares.}
\label{remonly}
\end{figure}

\begin{figure}[h]
\centering
\def\svgwidth{350pt}
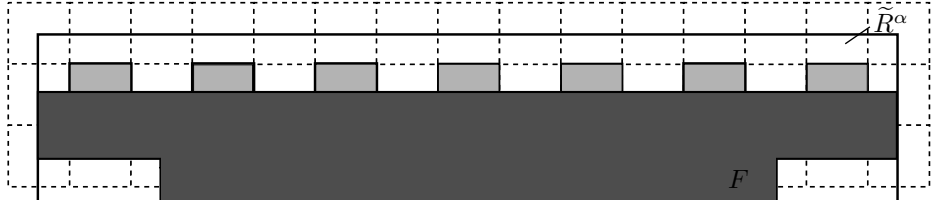
\caption{$F$ has a $\beta$-side.}
\label{albet}
\end{figure}

Another possibility is that $F$ has a $\beta$-type side, that is the portion of $\partial F$ intersecting the $\beta$-squares is a horizontal segment, as in Fig.~\ref{albet}.

\noindent
{\bf Step 3: conclusion.} We can repeat this procedure for each side, and finally, again by $\alpha$-rectangularization, we obtain that either $F$ is an $\alpha$-type rectangle or it has some $\beta$-type side.  In both cases, $F$ is a coordinate rectangle.
We note that all the estimates above can be iterated at each step and hold uniformly as long as the sides of $E_\ep^k$ are larger than $\eta$ (just to avoid that the length of any side vanishes), since they depend only on $c(\eta)$.
\endproof

\begin{figure}[ht]
\centering
\def\svgwidth{300pt}
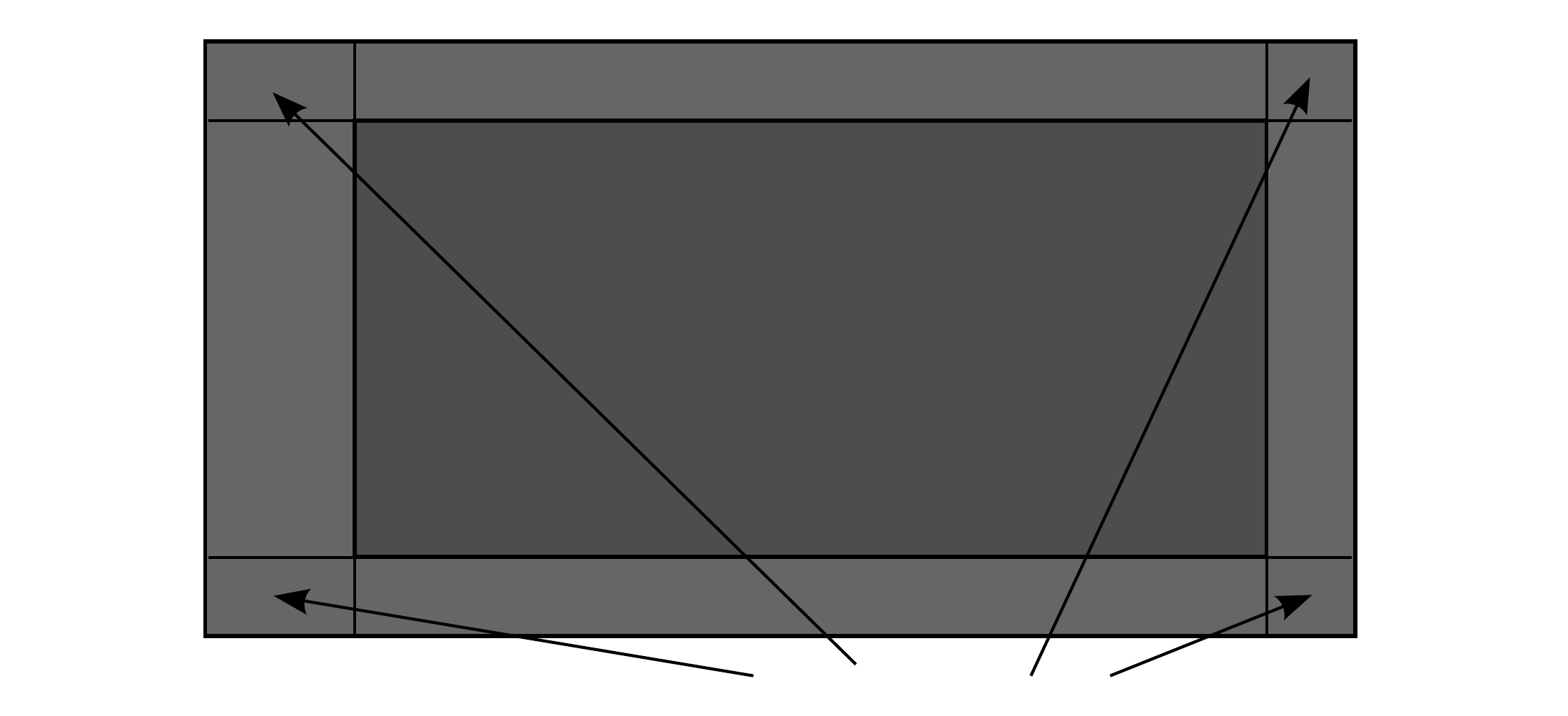
\caption{Picture of $E^{k+1}_\ep$ inside $E^k_\ep$}\label{BNG-figure4}
\end{figure}

As shown in \cite{BGN}, the motion of each side of $E_\ep^k$ can be studied separately remarking that the bulk term due to the small corner rectangles in Fig.~\ref{BNG-figure4} is negligible as $\ep\to0$.
As a consequence, we can describe the motion in terms of the length of the sides of $E_\ep^k$. This will be done in the following sections.

\subsection{A velocity function depending on $\delta$}\label{effveloc}

By the previous remark, we can reduce to a one-dimensional problem.
Let $x_k$ represents the projection of the left-hand vertical side of $E_k=E_\ep^k$ on the horizontal axis. The location of $x_{k+1}$ depends on a minimization argument involving $x_k$ and the length $L_k$ of the corresponding side of $E_k$. However, we will see that this latter dependence is locally constant, except for a discrete set of values of $L_k$. Indeed, for all $Y>0$ consider the minimum problems
\begin{equation}
\min\left\{g(N): N\in\mathbb{N}
\right\}
\label{minprob}
\end{equation}
where

\begin{equation}
g(N)=
\begin{cases}
-2\alpha N+\displaystyle\frac{N(N+1)}{2Y},&\mbox{$N$ even,}\\
\\
-2\alpha N+\displaystyle\frac{\delta\gamma}{Y}+\displaystyle\frac{N(N+1)}{2Y},& \mbox{$N$ odd}.
\end{cases}
\end{equation}
\\
Then the set of $Y>0$ for which (\ref{minprob}) does not have a unique solution is discrete. To check this it suffices to remark that the function to minimize is represented (up to multiplying by $2Y$) by two parabolas
\begin{equation*}
-4\alpha YX+X(X+1)\quad \text{ and }\quad-4\alpha YX+X(X+1)+2\delta\gamma
\end{equation*}
with minimum at
\begin{equation*}
X=2\alpha Y-\frac{1}{2}.
\end{equation*}
The minimizers in (\ref{minprob}) are not unique in the case that

\begin{equation}
g(N-1)=g(N)\quad\text{ or }\quad g(N)=g(N+1),
\end{equation}
\\
that is for $Y=\displaystyle\frac{N+\delta\gamma}{2\alpha}$ or $\displaystyle Y=\frac{N+1-\delta\gamma}{2\alpha}$ if $N$ is odd, while for $Y=\displaystyle\frac{N-\delta\gamma}{2\alpha}$ or $Y=\displaystyle\frac{N+1+\delta\gamma}{2\alpha}$ if $N$ is even.

\begin{definition}
We define the \emph{singular set $S_\delta$} for problems (\ref{minprob}) as
\begin{equation}
S_\delta=\frac{1}{2\alpha}\left[\left(2\mathbb{Z}+1+\delta\gamma\right)\cup\left(2\mathbb{Z}-\delta\gamma\right)\right].
\label{singular}
\end{equation}
\end{definition}

\begin{proposition}\label{thresh}
Let $Y\in(0,+\infty)\backslash S_\delta$ be fixed and $\widetilde{N}$ be the solution of the corresponding minimum problem (\ref{minprob}). Then there exists a threshold value of the contrast parameter
\begin{equation}
\widetilde{\delta}:=\frac{1}{2\gamma}
\label{thr}
\end{equation}
such that if $\delta\geq\widetilde{\delta}$  then $\widetilde{N}$ is even, while if $\delta<\widetilde{\delta}$ then $\widetilde{N}$ may be any integer.
\end{proposition}

\proof
Let $\widetilde{N}$ be odd. $\widetilde{N}$ is the unique solution in (\ref{minprob}), so that it satisfies
\begin{equation}
\begin{cases}
\displaystyle \widetilde{N}-\frac{1}{2}<2\alpha Y-\frac{1}{2}<\widetilde{N}+\frac{1}{2},\\
\\
g(\widetilde{N})<g(\widetilde{N}-1)\\
\\
g(\widetilde{N})<g(\widetilde{N}+1)
\end{cases}
\label{systembeta}
\end{equation}
that is,
\begin{equation*}
\begin{cases}
\displaystyle\frac{\widetilde{N}}{2\alpha}<Y<\frac{\widetilde{N}+1}{2\alpha}\\
\\
Y>\displaystyle\frac{\widetilde{N}+\delta\gamma}{2\alpha}\\
\\
Y<\displaystyle\frac{\widetilde{N}+1-\delta\gamma}{2\alpha}.
\end{cases}
\end{equation*}
We note that it is
\begin{equation*}
\frac{\widetilde{N}+\delta\gamma}{2\alpha}<\frac{\widetilde{N}+1-\delta\gamma}{2\alpha},
\end{equation*}
so that the system (\ref{systembeta}) has solutions, if and only if $\delta<\widetilde{\delta}$.
\endproof

Now we examine the iterated minimizing scheme for $\gamma/L_k=\gamma/L\in(0,+\infty)\backslash S_\delta$ fixed, which reads

\begin{equation}
\begin{cases}
x_{k+1}^L=x_k^L+\overline{N},\quad k\geq0\\
x_0^L=x^0
\end{cases}
\label{system}
\end{equation}
\\
with $x^0\in\{0,1\}$ and $\overline{N}\in\mathbb{N}$ the minimizer of

\begin{equation}
\min
\begin{cases}
-2\alpha N+\displaystyle\frac{1}{\gamma}\frac{N(N+1)}{2}L,&\text{$N$ even},\\
\\
-2\alpha N+\delta L+\displaystyle\frac{1}{\gamma}\frac{N(N+1)}{2}L,&\text{$N$ odd},
\end{cases}
\label{minimization}
\end{equation}
\\
which is unique up to the requirement that $\gamma/L\not\in S_\delta$.

\begin{remark}\label{perino}
As a trivial remark, after at most two steps $\{x_{k}^L\}_{k\geq0}$ is \emph{periodic modulo} $2$, that is, there exist integers $\bar{k},M\leq2$ and $n\geq1$ such that
\begin{equation}\label{periodi}
x_{k+M}^L=x_{k}^L+2n
\qquad\hbox{ for all } k\geq\bar{k}.
\end{equation}
For this, we note that $\{x_{k}^L\}_{k\geq0}$ is an arithmetic sequence and the conclusion depends whether $\overline{N}$ is odd or even. Moreover, the quotient $n/M$ depends only on $\gamma/L$ and $\delta$. In particular, if $\delta\geq1/2\gamma$ then $\bar{k}=M=1$.
By Proposition~\ref{thresh}, this is a straightforward consequence of Proposition~3.6 in \cite{BraSci13} with $N_\alpha=N_\beta=1$.
\end{remark}

\begin{definition}[effective velocity] \rm
We define the \emph{effective velocity function}\\ $f_\delta:(0,+\infty)\setminus S_\delta\longrightarrow[0,+\infty)$ by setting
\begin{equation}
f_\delta(Y)=\frac{2n}{M},
\label{velocityfunction}
\end{equation}
with $M$ and $n$ in (\ref{periodi}) defined by $L$ and $\gamma$ such that $Y=\gamma/L$. By Remark~\ref{perino} this is a good definition.
\label{effvel}
\end{definition}

We recall some properties of the velocity function (for the proof see \cite{BraSci13}).

\begin{remark}[properties of the velocity function $f_\delta$] The velocity function $f_\delta$ has the following properties:
\begin{description}
\item[(a)] $f_\delta$ is constant on each interval contained in its domain;
\item[(b)] $f_\delta(Y)=0$ if $$
Y<\overline{Y}_\delta:=\min\left\{\frac{3}{4\alpha},\frac{\delta\gamma+1}{2\alpha}\right\},
$$
where $\overline{Y}_\delta=\gamma/\overline{L}_\delta$ and $\overline{L}_\delta$ is the pinning threshold (see Subsection~\ref{newpinning}).\\ In particular,
$$
\lim_{\gamma\to 0^+} {\frac{1}{\gamma}} f_\delta\Bigl({\frac{\gamma}{L}}\Bigr)=0\,;
$$
\item[(c)] $f_\delta(Y)$ is an integer value;
\item[(d)] $f_\delta(Y)$ is a non decreasing function of $Y$;
\item[(e)] we have
$$
\lim_{\gamma\to+\infty} {\frac{1}{\gamma}} f_\delta\Bigl({\frac{\gamma}{L}}\Bigr)={\frac{2\alpha}{L}}\,.
$$
\end{description}
\end{remark}

\subsection{The effective pinning threshold}\label{newpinning}

We now examine the case when the limit motion is trivial; i.e., all $E_k = E^k_\ep$ are the same after a finite number of steps. This will be done by computing the \emph{pinning threshold}; i.e., the critical value of the side length $L$ above which it is energetically not favorable for a side to move.

\begin{figure}[h]
\centering
\def\svgwidth{200pt}
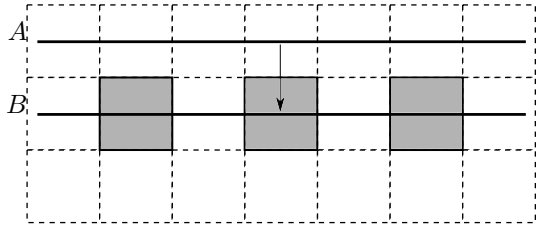
\caption{If $\delta<\widetilde{\delta}$ the motion is possible if the side can move at least by $\ep$.}
\label{pinning2}
\end{figure}

If $0\leq\delta<\widetilde{\delta}$ to compute it we have to impose that it is not energetically favorable to move inward a side by $\ep$. We then write the variation of the energy functional $\mathcal{F}^{\alpha,\beta_\ep}_{\ep,\tau}$ from configuration $A$ to configuration $B$ in Fig.~\ref{pinning2}, regarding a side of length $L$. If we impose it to be positive, we have
\begin{equation*}
-2\alpha\ep+L(\beta_\ep-\alpha)+\frac{1}{\tau}L\ep^2=\ep\left[-2\alpha+L\left(\delta+\frac{1}{\gamma}\right)\right]\geq0
\end{equation*}
and we obtain the pinning threshold
\begin{equation}
\widetilde{L}_\delta:=\frac{2\alpha\gamma}{\delta\gamma+1}.
\end{equation}
Note that if $\delta=0$ (i.e., $\beta_\ep=\alpha$), then we recover the threshold of the homogeneous case
\begin{equation*}
\widetilde{L}_0=2\alpha\gamma.
\end{equation*}

\begin{figure}[h]
\centering
\def\svgwidth{200pt}
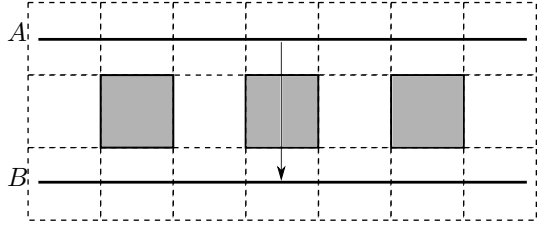
\caption{If $\delta\geq \widetilde{\delta}$ the motion is possible if the side can move at least by $2\ep$.}
\label{pinning3}
\end{figure}

If $\delta\geq\widetilde{\delta}$, instead, by the condition that $E_k$ be an $\alpha$-type rectangle, we have to impose that it is not energetically favorable to move inward a side by $2\ep$ (see Fig.~\ref{pinning3}). As shown in \cite{BraSci13}, this gives the pinning threshold
\begin{equation*}
\widetilde{L}_{\widetilde{\delta}}=\frac{4}{3}\alpha\gamma.
\end{equation*}

\begin{figure}[h]
\centering
\includegraphics[scale=.90]{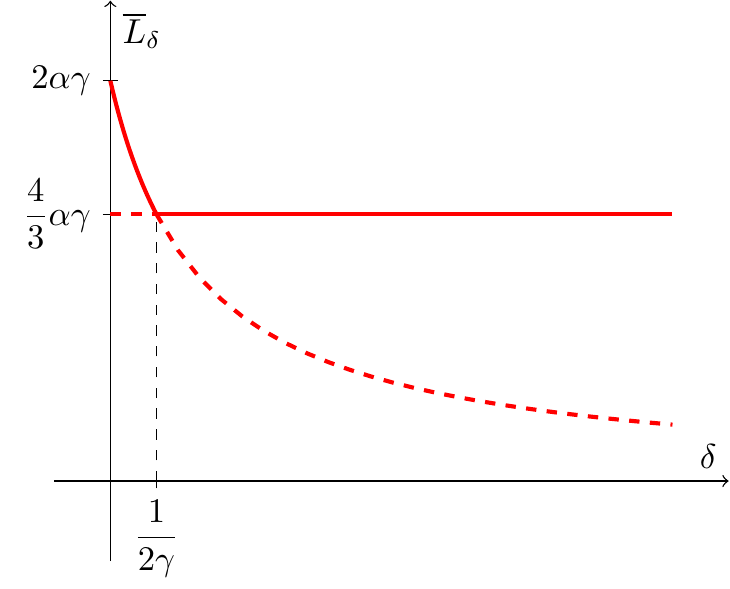}
\caption{Effective pinning threshold.}
\label{threshold}
\end{figure}

Hence, the \emph{effective pinning threshold} (see Fig.~\ref{threshold}) is given by
\begin{equation}
\overline{L}_\delta=\max\left\{\frac{2\alpha\gamma}{\delta\gamma+1},\frac{4}{3}\alpha\gamma\right\}.
\label{pinningeff}
\end{equation}

\subsection{Computation of the velocity function}\label{computation}
In this section we compute explicitly the velocity function $f_\delta$ assuming, without loss of generality, that $\gamma=1$. We restrict ourselves to the case $\delta<1/2$, because if $\delta\geq1/2$ the velocity function is given by (see \cite{BraSci13}, Section~4)
\begin{equation}
\overline{f}(Y)=2\left\lfloor \alpha Y+\frac{1}{4}\right\rfloor.
\label{highvel}
\end{equation}
We denote by $\overline{N}$ the minimizer of the problem (\ref{minimization}) and subdivide the computation into different cases:
\\
\\
(a) $x_n$ is even and $Y\in\left(\displaystyle\frac{2k+1+\delta}{2\alpha},\frac{2k+2-\delta}{2\alpha}\right)$ for some $k\geq0$; in this case $\overline{N}=2k+1$ and $x_{n+1}= x_n+\overline{N}$ is odd. The next point is $x_{n+2}=x_{n+1}+\overline{N}=x_n+2\overline{N}$, which is even, so that the sequence $\{x_m\}$ oscillates between even and odd numbers (that is, the side is alternatively $\alpha$-type  and $\beta$-type). In this case,

\begin{equation*}
f_\delta(Y)=\frac{x_{n+2}-x_n}{2}=\frac{2\overline{N}}{2}=2k+1=\lfloor 2\alpha Y\rfloor;
\end{equation*}
\\
(b) $x_n$ is odd and $Y\in\left(\displaystyle\frac{2k+1+\delta}{2\alpha},\frac{2k+2-\delta}{2\alpha}\right)$ for some $k\geq0$; in this case
$x_{n+1}=x_n+\overline{N}$ is even and $x_{n+2}=x_n+2\overline{N}$, is odd, so that as before

\begin{equation*}
f_\delta(Y)=2k+1=\lfloor 2\alpha Y\rfloor;
\end{equation*}
\\
(c) $x_n$ is even and $Y\in\left(\displaystyle\frac{2k-\delta}{2\alpha},\frac{2k+1+\delta}{2\alpha}\right)$ for some $k\geq0$; in this case $\overline{N}=2k$ and $x_{n+1}= x_n+\overline{N}$ is even. Therefore the sequence $\{x_m\}$ consists of only even numbers (that is, at each step the side is $\alpha$-type) and in this case the velocity function is given by

\begin{equation*}
f_\delta(Y)=x_{n+1}-x_n=\overline{N}=2k;
\end{equation*}
\\
(d) $x_n$ is odd and $Y\in\left(\displaystyle\frac{2k-\delta}{2\alpha},\frac{2k+1+\delta}{2\alpha}\right)$ for some $k\geq0$; in this case
$x_{n+1}= x_n+\overline{N}$ is also odd. Therefore the sequence $\{x_m\}$ consists of only odd numbers (that is, at each step the side is $\beta$-type) and in this case the velocity function is  given again by

\begin{equation*}
f_\delta(Y)=x_{n+1}-x_n=N=2k.
\end{equation*}
\\
Note that, collecting all the cases, we can write the velocity function as
\begin{equation}
f_\delta(Y)=
\begin{cases}
0&\text{if }0<Y<\displaystyle\frac{\delta+1}{2\alpha},\\
\\
2k& \text{if }Y\in\left(\displaystyle\frac{2k-\delta}{2\alpha},\frac{2k+1+\delta}{2\alpha}\right),\qquad k\geq0\\
\\
2k+1& \text{if }Y\in\left(\displaystyle\frac{2k+1+\delta}{2\alpha},\frac{2k+2-\delta}{2\alpha}\right).
\end{cases}
\label{fdelta}
\end{equation}
It can be rewritten equivalently as
\begin{equation*}
f_\delta(Y)=
\begin{cases}
\lfloor 2\alpha Y\rfloor+1& \text{if }Y\in\left(\displaystyle\frac{2k-\delta}{2\alpha},\frac{2k}{2\alpha}\right),\\
\\
\lfloor 2\alpha Y\rfloor& \text{if }Y\in\left(\displaystyle\frac{2k}{2\alpha},\frac{2k+1}{2\alpha}\right),\\
\\
\lfloor 2\alpha Y\rfloor-1& \text{if }Y\in\left(\displaystyle\frac{2k+1}{2\alpha},\frac{2k+1+\delta}{2\alpha}\right),\\
\\
\lfloor 2\alpha Y\rfloor& \text{if }Y\in\left(\displaystyle\frac{2k+1+\delta}{2\alpha},\frac{2k+2-\delta}{2\alpha}\right).
\end{cases}
\end{equation*}
Therefore we notice accelerating and decelerating effects (due to the microstructure through $\delta$) with respect to the velocity function $\widetilde{f}$ obtained in the homogeneous case~\cite{BGN}, that is
\begin{equation*}
\widetilde{f}(Y)=
\begin{cases}
0&\text{if }0<Y<\displaystyle\frac{1}{2\alpha},\\
\\
\lfloor 2\alpha Y\rfloor &\text{if }Y\in\left(\displaystyle\frac{k}{2\alpha},\frac{k+1}{2\alpha}\right), k\geq1.
\end{cases}
\end{equation*}
Moreover, we recover $\widetilde{f}$ computing $f_\delta$ for $\delta=0$. If we choose $\delta=1/2$ (actually, any $\delta\geq1/2$), we recover the velocity function $\overline{f}$ (\ref{highvel}) which corresponds to the high\hbox{-}contrast case.\\

We conclude this section by writing the general formula of the velocity function $f_\delta$ valid for any $\delta$ and $\gamma$:

\begin{equation*}
f_\delta(Y)=
\begin{cases}
0&\text{if }0<Y<\displaystyle\frac{\min\{\delta\gamma,1/2\}+1}{2\alpha},\\
\\
2k& \text{if }Y\in\left(\displaystyle\frac{2k-\min\{\delta\gamma,1/2\}}{2\alpha},\frac{2k+1+\min\{\delta\gamma,1/2\}}{2\alpha}\right),\\
\\
2k+1& \text{if }Y\in\left(\displaystyle\frac{2k+1+\min\{\delta\gamma,1/2\}}{2\alpha},\frac{2k+2-\min\{\delta\gamma,1/2\}}{2\alpha}\right),
\end{cases}
\label{fdelta2}
\end{equation*}
with $k\geq0$.

\subsection{Description of the homogenized limit motion}\label{limitmotion}

The following characterization of any limit motion holds (see Theorem 3.11 in \cite{BraSci13}).

\begin{theorem}\label{limitmotion1}For all $\ep>0$, let $E^0_\ep\in\mathcal{D}_\ep$ be a coordinate rectangle with sides $S^0_{1,\ep},\dots,S^0_{4,\ep}$. Assume also that
\begin{equation*}
\lim_{\ep\to0^+}\emph{d}_\mathcal{H}(E^0_\ep,E_0)=0
\end{equation*}
for some fixed coordinate rectangle $E_0$. Let $\delta,\gamma>0$ be fixed and let  $E_\ep(t)= E_{\ep,\gamma\ep}(t)$ be the piecewise-constant motion with initial datum $E^0_\ep$ defined in {\rm(\ref{disefo})}.
Then, up to a subsequence, $E_\ep(t)$ converges as $\ep\to0$ to $E(t)$, where $E(t)$ is a coordinate rectangle with sides $S_i(t)$ and such that $E(0)=E_0$. Each $S_i$ moves inward with velocity $v_i(t)$ satisfying
\begin{equation}
v_i(t)\in\left[\displaystyle\frac{1}{\gamma}f_\delta\biggl({\frac{\gamma}{L_i(t)}}\biggr)^-,\displaystyle\frac{1}{\gamma}f_\delta\biggl({\frac{\gamma}{L_i(t)}}\biggr)^+\right],
\label{vl}
\end{equation}
where $f_\delta$ is given by Definition {\rm\ref{effvel}}, $L_i(t):=\mathcal{H}^1(S_i(t))$ denotes the length of the side $S_i(t)$, until the extinction time when $L_i(t)=0$, and $f_\delta(Y)^-,f_\delta(Y)^+$ are the lower and upper limits of the effective-velocity function at $Y\in (0,+\infty)$.
\end{theorem}

In case of a unique evolution, the limit motion is described as follows (see Theorem 3.12 in \cite{BraSci13}).

\begin{theorem}[unique limit motion]
Let $E_\ep,E_0$ be as in the statement of Theorem {\rm \ref{limitmotion1}}. Assume in addition that the lengths $L^0_1,L^0_2$ of the sides of the initial set $E_0$ satisfy one of the following conditions (we assume that $L^0_1\leq L^0_2$):

\begin{itemize}
\item[\emph{(a)}] $L^0_1,L^0_2>\overline{L}_\delta$, $\overline{L}_\delta$ given by {\rm(\ref{pinningeff})} \emph{(total pinning)};
\item[\emph{(b)}] $L^0_1<\overline{L}_\delta$ and $L^0_2\leq\overline{L}_\delta$ \emph{(vanishing in finite time)};
\end{itemize}
then $E_\ep(t)$ converges locally in time to $E(t)$ as $\ep\to 0$, where $E(t)$ is the unique rectangle with sides of lengths $L_1(t)$ and $L_2(t)$ which solve the following system of ordinary differential equations
\begin{equation}\label{unita}
\begin{cases}\displaystyle
\dot{L}_1(t)=-{\frac{2}{\gamma}}\,f_\delta\left({\frac{\gamma}{L_2(t)}}\right)\\
\\ \displaystyle
\dot{L}_2(t)=-{\frac{2}{\gamma}}\, f_\delta\left({\frac{\gamma}{L_1(t)}}\right)
\end{cases}
\end{equation}
for almost every $t$, with initial conditions $L_1(0)=L^0_1$ and $L_2(0)=L^0_2$.
\end{theorem}

\section{The periodic case with $K$ contrast parameters}\label{periodic2}
In this section we study the same problem as before in a more general framework. We consider a medium with inclusions distributed into periodic uniform layers as follows.\\
\begin{figure}[h]
\centering
\def\svgwidth{150pt}
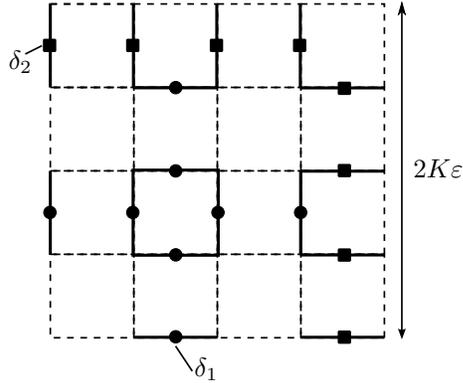
\caption{The periodicity cell for $K=2$.}
\label{cella2}
\end{figure}
Let $\ep>0$ be fixed and $\delta_1,\delta_2,\dots,\delta_K,K\in\mathbb{N}$ be positive. We consider $2K\ep$-periodic coefficients $c_{ij}^\ep$ indexed on nearest-neighbors of $\ep\mathbb{Z}^2$ and defined for $i,j$ such that
\begin{equation*}
0\leq \frac{i_1+j_1}{2},\frac{i_2+j_2}{2}<2K\ep
\end{equation*}
by
\begin{equation}
c_{ij}^\ep=
\begin{cases}
\alpha+\delta_r\ep, & \text{if }\displaystyle\frac{i_1+j_1}{2},\frac{i_2+j_2}{2}=\left(2r-\frac{1}{2}\right)\ep,\quad r=1,\dots,K\\
\\
\alpha,   & \text{otherwise.}
\end{cases}
\end{equation}
In Fig.~\ref{cella2} the periodicity cell is pictured in the case $K=2$. Here the bonds with parameter $\delta_1$ are marked with a dot, the ones with parameter $\delta_2$ are marked with a square and the dashed lines represent the $\alpha$-bonds.

Correspondingly, to these coefficients we associate the energy $\text{P}^{\alpha,\beta_\ep}_\ep(\mathcal{I})$
defined on subsets $\mathcal{I}$ of $\ep\mathbb{Z}^2$ as in (\ref{energy}). We consider the same discrete-in-time minimization scheme for the energy $\mathcal{F}_{\ep,\tau}^{\alpha,\beta_\ep}$ with $\tau=\gamma\ep$ as in Subsection~\ref{timemin} and we restrict our analysis to rectangular evolutions as in Section~\ref{rectangle}. We will see that the minimization problem and the velocity function depend on the choice of $\delta_r,r=1,\dots,K$; in particular, on their relative position with respect to the critical value $\widetilde{\delta}$ defined by equation (\ref{thr}).

We will treat only the cases
\begin{equation}
\widetilde{\delta}\leq\delta_r\quad\text{for some $r\in\{1,\dots,K\}$}
\label{ass1}
\end{equation}
and
\begin{equation}
0\leq\delta_r<\widetilde{\delta}, \quad \forall r=1,\dots,K,
\label{ass2}
\end{equation}
because if $\widetilde{\delta}\leq\delta_r$ for all $r$ then we are in the high-contrast case already described in \cite{BraSci13}.

\subsection{The pinning threshold}
For the computation of the pinning threshold we refer to Subsection~\ref{newpinning}.

Under assumption (\ref{ass1}), after a finite number of steps the side is pinned if it cannot move inward by $2\ep$. In this case, the pinning threshold is given by
\begin{equation*}
\widetilde{L}_{\widetilde{\delta}}=\frac{4}{3}\alpha\gamma.
\end{equation*}

If (\ref{ass2}) holds, instead, after a finite number of steps the side is pinned if it cannot move inward by $\ep$. In particular, the pinning threshold now depends on $\delta_{\bar{r}}=\displaystyle\min_{1\leq r\leq K}\{\delta_r\}$ and it is given by
\begin{equation*}
\widetilde{L}_{\delta_{\bar{r}}}=\frac{2\alpha\gamma}{\delta_{\bar{r}}\gamma+1}.
\end{equation*}
Hence, collecting the two cases we obtain the pinning threshold
\begin{equation}
\overline{L}_{\delta_1,\dots,\delta_K}=\max\{\widetilde{L}_{\delta_{\bar{r}}},\widetilde{L}_{\widetilde{\delta}}\}.
\end{equation}

\subsection{The effective velocity function}
We adopt the same notation as in Subsection \ref{effveloc}. For all $Y>0$ we consider the minimum problems

\begin{equation}
\min\left\{g(N): N\in\mathbb{N}\right\}
\label{minprob1}
\end{equation}
\\
where

\begin{equation}
g(N)=
\begin{cases}
-2\alpha N+\displaystyle\frac{N(N+1)}{2Y},& \text{if }[N]_{2K}=[2r-2]_{2K},\\
\\
-2\alpha N+\displaystyle\frac{\delta_r\gamma}{Y}+\displaystyle\frac{N(N+1)}{2Y},&\text{if } [N]_{2K}=[2r-1]_{2K},r=1,\dots,K,
\end{cases}
\end{equation}
\\
and $[z]_{2K}$ is the congruence class of $z$ modulo $2K$.
Then the set of $Y>0$ for which (\ref{minprob1}) does not have a unique solution is discrete. For this we remark that the function to minimize is represented by $K+1$ parabolas
\begin{equation*}
-4\alpha YX+X(X+1)\quad\text{ and }\quad-4\alpha YX+X(X+1)+2\delta_r\gamma\quad r=1,\dots,K
\end{equation*}
with {minimum at}
\begin{equation*}
X=2\alpha Y-\frac{1}{2}.
\end{equation*}
As a consequence of (\ref{singular}) we have that the minimizers in (\ref{minprob1}) are not unique in the case that $Y\in S_{\delta_r}, r=1,\dots,K$ where
\begin{equation}
S_{\delta_r}:={\frac{1}{2\alpha}}\left[(2(2r-1)K\mathbb{Z}+C_{\delta_r\gamma})\cup((2(2r-1)K\mathbb{Z}+1-C_{\delta_r\gamma})\right]
\label{sing}
\end{equation}
and $C_{\delta_r\gamma}=\min\{\delta_r\gamma,1/2\}, r=1,\dots,K$.

\begin{definition}
We define the \emph{singular set {$S_{\delta_1,\dots,\delta_K}$}} for problems (\ref{minprob1}) as
\begin{equation}
S_{\delta_1,\dots,\delta_K}=\displaystyle\bigcup_{r=1}^K S_{\delta_r}
\end{equation}where $S_{\delta_r}$ is defined by (\ref{sing}).
\label{singular2}
\end{definition}

We now examine the iterated minimizing scheme for $\gamma/L\in(0,+\infty)\backslash S_{\delta_1,\dots,\delta_K}$ fixed, which reads
\begin{equation}
\begin{cases}
x_{k+1}^L=x_k^L+\overline{N},\quad k\geq0\\
x_0^L=x^0
\end{cases}
\label{system2}
\end{equation}
with $x^0\in\{0,1,2,\dots,2K-1\}$ and $\overline{N}\in\mathbb{N}$ the minimizer of
\begin{equation*}
\min
\begin{cases}
-2\alpha N+\displaystyle\frac{1}{\gamma}\frac{N(N+1)}{2}L,& \text{if }[N]_{2K}=[2r-2]_{2K},\\
\\
-2\alpha N+\delta_r L+\displaystyle\frac{1}{\gamma}\frac{N(N+1)}{2}L,& \text{if }[N]_{2K}=[2r-1]_{2K},r=1,\dots,K,
\end{cases}
\label{minimization2}
\end{equation*}
which is unique up to the requirement that $\gamma/L\not\in S_{\delta_1,\dots,\delta_K}$. With an analogous argument as in Section~\ref{effveloc} we can prove that, after at most $2K$ steps, $\{x_{k}^L\}_{k\geq0}$ is periodic modulo $2K$. Hence, we can define the effective velocity function $f=f_{\delta_1,\dots,\delta_K}$ as in Definition~\ref{effvel}.

\subsection{Computation of the velocity function}
In this section we give the expression of the velocity function without proof, which follows by analogous computations as in Subsection~\ref{computation}.

For any $\gamma,\delta_1,\dots,\delta_K$, the velocity function $f=f_{\delta_1,\dots,\delta_K}$ is given by

\begin{equation*}
f(Y)=
\begin{cases}
0, & \text{if }0<Y<\gamma/\overline{L}_{\delta_1,\dots,\delta_K},\\
\\
2Kk, & \text{if }Y\in\left(\displaystyle\frac{2Kk-C_{\delta_K\gamma}}{2\alpha},\frac{2Kk+1+C_{\delta_1\gamma}}{2\alpha}\right),\\
\\
2Kk+1, & \text{if }Y\in\left(\displaystyle\frac{2Kk+1+C_{\delta_1\gamma}}{2\alpha},\frac{2Kk+2-C_{\delta_1\gamma}}{2\alpha}\right),\\
\\
2Kk+2, & \text{if }Y\in\left(\displaystyle\frac{2Kk+2-C_{\delta_1\gamma}}{2\alpha},\frac{2Kk+3+C_{\delta_2\gamma}}{2\alpha}\right),\\
\\
\quad\vdots&\quad\vdots\qquad\qquad\vdots\\
\\
2Kk+2K-1, & \text{if }Y\in\left(\displaystyle\frac{2K(k+1)-1+C_{\delta_K\gamma}}{2\alpha},\frac{2K(k+1)-C_{\delta_K\gamma}}{2\alpha}\right),
\end{cases}
\end{equation*}
with $k\geq0$.

\section*{Acknowledgments} I am grateful to Andrea Braides for suggesting this problem, and I would like to thank him for his advices. I acknowledge the anonymous referees for their interesting remarks leading to improvements of the manuscript.



\end{document}

%% file: fig0.pdf_tex
\begingroup%
  \makeatletter%
  \providecommand\color[2][]{%
    \errmessage{(Inkscape) Color is used for the text in Inkscape, but the package 'color.sty' is not loaded}%
    \renewcommand\color[2][]{}%
  }%
  \providecommand\transparent[1]{%
    \errmessage{(Inkscape) Transparency is used (non-zero) for the text in Inkscape, but the package 'transparent.sty' is not loaded}%
    \renewcommand\transparent[1]{}%
  }%
  \providecommand\rotatebox[2]{#2}%
  \ifx\svgwidth\undefined%
    \setlength{\unitlength}{151.96165771bp}%
    \ifx\svgscale\undefined%
      \relax%
    \else%
      \setlength{\unitlength}{\unitlength * \real{\svgscale}}%
    \fi%
  \else%
    \setlength{\unitlength}{\svgwidth}%
  \fi%
  \global\let\svgwidth\undefined%
  \global\let\svgscale\undefined%
  \makeatother%
  \begin{picture}(1,1.16349347)%
    \put(0,0){\includegraphics[width=\unitlength]{fig0.pdf}}%
    \put(0.24896788,0.04190544){\makebox(0,0)[lb]{\smash{$N_\beta$}}}%
    \put(0.74856547,0.04190544){\makebox(0,0)[lb]{\smash{$N_\alpha$}}}%
  \end{picture}%
\endgroup%

%% file: heuristic.pdf_tex
\begingroup%
  \makeatletter%
  \providecommand\color[2][]{%
    \errmessage{(Inkscape) Color is used for the text in Inkscape, but the package 'color.sty' is not loaded}%
    \renewcommand\color[2][]{}%
  }%
  \providecommand\transparent[1]{%
    \errmessage{(Inkscape) Transparency is used (non-zero) for the text in Inkscape, but the package 'transparent.sty' is not loaded}%
    \renewcommand\transparent[1]{}%
  }%
  \providecommand\rotatebox[2]{#2}%
  \ifx\svgwidth\undefined%
    \setlength{\unitlength}{308.8bp}%
    \ifx\svgscale\undefined%
      \relax%
    \else%
      \setlength{\unitlength}{\unitlength * \real{\svgscale}}%
    \fi%
  \else%
    \setlength{\unitlength}{\svgwidth}%
  \fi%
  \global\let\svgwidth\undefined%
  \global\let\svgscale\undefined%
  \makeatother%
  \begin{picture}(1,0.43305052)%
    \put(0,0){\includegraphics[width=\unitlength]{heuristic.pdf}}%
    \put(0.47861411,0.39223405){\color[rgb]{0,0,0}\makebox(0,0)[lb]{\smash{$L$}}}%
  \end{picture}%
\endgroup%

%% file: figra1.pdf_tex
\begingroup%
  \makeatletter%
  \providecommand\color[2][]{%
    \errmessage{(Inkscape) Color is used for the text in Inkscape, but the package 'color.sty' is not loaded}%
    \renewcommand\color[2][]{}%
  }%
  \providecommand\transparent[1]{%
    \errmessage{(Inkscape) Transparency is used (non-zero) for the text in Inkscape, but the package 'transparent.sty' is not loaded}%
    \renewcommand\transparent[1]{}%
  }%
  \providecommand\rotatebox[2]{#2}%
  \ifx\svgwidth\undefined%
    \setlength{\unitlength}{105.78056641bp}%
    \ifx\svgscale\undefined%
      \relax%
    \else%
      \setlength{\unitlength}{\unitlength * \real{\svgscale}}%
    \fi%
  \else%
    \setlength{\unitlength}{\svgwidth}%
  \fi%
  \global\let\svgwidth\undefined%
  \global\let\svgscale\undefined%
  \makeatother%
  \begin{picture}(1,0.96736223)%
    \put(0,0){\includegraphics[width=\unitlength]{figra1.pdf}}%
    \put(0.2361378,-0.03002088){\color[rgb]{0,0,0}\makebox(0,0)[lb]{\smash{$\beta_\varepsilon$}}}%
    \put(0.71713299,-0.03002088){\color[rgb]{0,0,0}\makebox(0,0)[lb]{\smash{$\alpha$}}}%
    \put(-0.00302572,0.26469888){\color[rgb]{0,0,0}\makebox(0,0)[lb]{\smash{$\varepsilon$}}}%
  \end{picture}%
\endgroup%

%% file: remove.pdf_tex
\begingroup%
  \makeatletter%
  \providecommand\color[2][]{%
    \errmessage{(Inkscape) Color is used for the text in Inkscape, but the package 'color.sty' is not loaded}%
    \renewcommand\color[2][]{}%
  }%
  \providecommand\transparent[1]{%
    \errmessage{(Inkscape) Transparency is used (non-zero) for the text in Inkscape, but the package 'transparent.sty' is not loaded}%
    \renewcommand\transparent[1]{}%
  }%
  \providecommand\rotatebox[2]{#2}%
  \ifx\svgwidth\undefined%
    \setlength{\unitlength}{333.03989258bp}%
    \ifx\svgscale\undefined%
      \relax%
    \else%
      \setlength{\unitlength}{\unitlength * \real{\svgscale}}%
    \fi%
  \else%
    \setlength{\unitlength}{\svgwidth}%
  \fi%
  \global\let\svgwidth\undefined%
  \global\let\svgscale\undefined%
  \makeatother%
  \begin{picture}(1,0.41206474)%
    \put(0,0){\includegraphics[width=\unitlength]{remove.pdf}}%
    \put(0.31999571,0.23320235){\color[rgb]{0,0,0}\makebox(0,0)[lb]{\smash{$\Gamma$}}}%
    \put(0.40338149,0.02104083){\color[rgb]{0,0,0}\makebox(0,0)[lb]{\smash{$\bf\emph{F}$}}}%
    \put(0.94673965,0.02104083){\color[rgb]{0,0,0}\makebox(0,0)[lb]{\smash{$\bf\emph{F}$}}}%
  \end{picture}%
\endgroup%

%% file: add.pdf_tex
\begingroup%
  \makeatletter%
  \providecommand\color[2][]{%
    \errmessage{(Inkscape) Color is used for the text in Inkscape, but the package 'color.sty' is not loaded}%
    \renewcommand\color[2][]{}%
  }%
  \providecommand\transparent[1]{%
    \errmessage{(Inkscape) Transparency is used (non-zero) for the text in Inkscape, but the package 'transparent.sty' is not loaded}%
    \renewcommand\transparent[1]{}%
  }%
  \providecommand\rotatebox[2]{#2}%
  \ifx\svgwidth\undefined%
    \setlength{\unitlength}{333.03989258bp}%
    \ifx\svgscale\undefined%
      \relax%
    \else%
      \setlength{\unitlength}{\unitlength * \real{\svgscale}}%
    \fi%
  \else%
    \setlength{\unitlength}{\svgwidth}%
  \fi%
  \global\let\svgwidth\undefined%
  \global\let\svgscale\undefined%
  \makeatother%
  \begin{picture}(1,0.41206474)%
    \put(0,0){\includegraphics[width=\unitlength]{add.pdf}}%
    \put(0.31999571,0.23320235){\color[rgb]{0,0,0}\makebox(0,0)[lb]{\smash{$\Gamma$}}}%
    \put(0.40338149,0.02104083){\color[rgb]{0,0,0}\makebox(0,0)[lb]{\smash{$\bf\emph{F}$}}}%
    \put(0.94673965,0.02104083){\color[rgb]{0,0,0}\makebox(0,0)[lb]{\smash{$\bf\emph{F}$}}}%
  \end{picture}%
\endgroup%

%% file: maximal.pdf_tex
\begingroup%
  \makeatletter%
  \providecommand\color[2][]{%
    \errmessage{(Inkscape) Color is used for the text in Inkscape, but the package 'color.sty' is not loaded}%
    \renewcommand\color[2][]{}%
  }%
  \providecommand\transparent[1]{%
    \errmessage{(Inkscape) Transparency is used (non-zero) for the text in Inkscape, but the package 'transparent.sty' is not loaded}%
    \renewcommand\transparent[1]{}%
  }%
  \providecommand\rotatebox[2]{#2}%
  \ifx\svgwidth\undefined%
    \setlength{\unitlength}{450.07998047bp}%
    \ifx\svgscale\undefined%
      \relax%
    \else%
      \setlength{\unitlength}{\unitlength * \real{\svgscale}}%
    \fi%
  \else%
    \setlength{\unitlength}{\svgwidth}%
  \fi%
  \global\let\svgwidth\undefined%
  \global\let\svgscale\undefined%
  \makeatother%
  \begin{picture}(1,0.21418415)%
    \put(0,0){\includegraphics[width=\unitlength]{maximal.pdf}}%
    \put(0.938604,0.18133746){\color[rgb]{0,0,0}\makebox(0,0)[lb]{\smash{$\widetilde{R}^\alpha$}}}%
    \put(0.78105986,0.02270522){\color[rgb]{0,0,0}\makebox(0,0)[lb]{\smash{$F$}}}%
  \end{picture}%
\endgroup%

%% file: repl.pdf_tex
\begingroup%
  \makeatletter%
  \providecommand\color[2][]{%
    \errmessage{(Inkscape) Color is used for the text in Inkscape, but the package 'color.sty' is not loaded}%
    \renewcommand\color[2][]{}%
  }%
  \providecommand\transparent[1]{%
    \errmessage{(Inkscape) Transparency is used (non-zero) for the text in Inkscape, but the package 'transparent.sty' is not loaded}%
    \renewcommand\transparent[1]{}%
  }%
  \providecommand\rotatebox[2]{#2}%
  \ifx\svgwidth\undefined%
    \setlength{\unitlength}{450.07998047bp}%
    \ifx\svgscale\undefined%
      \relax%
    \else%
      \setlength{\unitlength}{\unitlength * \real{\svgscale}}%
    \fi%
  \else%
    \setlength{\unitlength}{\svgwidth}%
  \fi%
  \global\let\svgwidth\undefined%
  \global\let\svgscale\undefined%
  \makeatother%
  \begin{picture}(1,0.21418415)%
    \put(0,0){\includegraphics[width=\unitlength]{repl.pdf}}%
    \put(0.938604,0.18133746){\color[rgb]{0,0,0}\makebox(0,0)[lb]{\smash{$\widetilde{R}^\alpha$}}}%
    \put(0.66788891,0.14916335){\color[rgb]{0,0,0}\makebox(0,0)[lb]{\smash{$R$}}}%
    \put(0.78105986,0.02270522){\color[rgb]{0,0,0}\makebox(0,0)[lb]{\smash{$F$}}}%
  \end{picture}%
\endgroup%

%% file: final.pdf_tex
\begingroup%
  \makeatletter%
  \providecommand\color[2][]{%
    \errmessage{(Inkscape) Color is used for the text in Inkscape, but the package 'color.sty' is not loaded}%
    \renewcommand\color[2][]{}%
  }%
  \providecommand\transparent[1]{%
    \errmessage{(Inkscape) Transparency is used (non-zero) for the text in Inkscape, but the package 'transparent.sty' is not loaded}%
    \renewcommand\transparent[1]{}%
  }%
  \providecommand\rotatebox[2]{#2}%
  \ifx\svgwidth\undefined%
    \setlength{\unitlength}{450.07998047bp}%
    \ifx\svgscale\undefined%
      \relax%
    \else%
      \setlength{\unitlength}{\unitlength * \real{\svgscale}}%
    \fi%
  \else%
    \setlength{\unitlength}{\svgwidth}%
  \fi%
  \global\let\svgwidth\undefined%
  \global\let\svgscale\undefined%
  \makeatother%
  \begin{picture}(1,0.21418415)%
    \put(0,0){\includegraphics[width=\unitlength]{final.pdf}}%
    \put(0.938604,0.18133746){\color[rgb]{0,0,0}\makebox(0,0)[lb]{\smash{$\widetilde{R}^\alpha$}}}%
    \put(0.78208323,0.04443655){\color[rgb]{0,0,0}\makebox(0,0)[lb]{\smash{$F$}}}%
  \end{picture}%
\endgroup%

%% file: onlyone.pdf_tex
\begingroup%
  \makeatletter%
  \providecommand\color[2][]{%
    \errmessage{(Inkscape) Color is used for the text in Inkscape, but the package 'color.sty' is not loaded}%
    \renewcommand\color[2][]{}%
  }%
  \providecommand\transparent[1]{%
    \errmessage{(Inkscape) Transparency is used (non-zero) for the text in Inkscape, but the package 'transparent.sty' is not loaded}%
    \renewcommand\transparent[1]{}%
  }%
  \providecommand\rotatebox[2]{#2}%
  \ifx\svgwidth\undefined%
    \setlength{\unitlength}{450.07998047bp}%
    \ifx\svgscale\undefined%
      \relax%
    \else%
      \setlength{\unitlength}{\unitlength * \real{\svgscale}}%
    \fi%
  \else%
    \setlength{\unitlength}{\svgwidth}%
  \fi%
  \global\let\svgwidth\undefined%
  \global\let\svgscale\undefined%
  \makeatother%
  \begin{picture}(1,0.21418415)%
    \put(0,0){\includegraphics[width=\unitlength]{onlyone.pdf}}%
    \put(0.938604,0.18133746){\color[rgb]{0,0,0}\makebox(0,0)[lb]{\smash{$\widetilde{R}^\alpha$}}}%
    \put(0.78105986,0.02270522){\color[rgb]{0,0,0}\makebox(0,0)[lb]{\smash{$F$}}}%
  \end{picture}%
\endgroup%

%% file: remonlyone.pdf_tex
\begingroup%
  \makeatletter%
  \providecommand\color[2][]{%
    \errmessage{(Inkscape) Color is used for the text in Inkscape, but the package 'color.sty' is not loaded}%
    \renewcommand\color[2][]{}%
  }%
  \providecommand\transparent[1]{%
    \errmessage{(Inkscape) Transparency is used (non-zero) for the text in Inkscape, but the package 'transparent.sty' is not loaded}%
    \renewcommand\transparent[1]{}%
  }%
  \providecommand\rotatebox[2]{#2}%
  \ifx\svgwidth\undefined%
    \setlength{\unitlength}{450.07998047bp}%
    \ifx\svgscale\undefined%
      \relax%
    \else%
      \setlength{\unitlength}{\unitlength * \real{\svgscale}}%
    \fi%
  \else%
    \setlength{\unitlength}{\svgwidth}%
  \fi%
  \global\let\svgwidth\undefined%
  \global\let\svgscale\undefined%
  \makeatother%
  \begin{picture}(1,0.21418415)%
    \put(0,0){\includegraphics[width=\unitlength]{remonlyone.pdf}}%
    \put(0.938604,0.18133746){\color[rgb]{0,0,0}\makebox(0,0)[lb]{\smash{$\widetilde{R}^\alpha$}}}%
    \put(0.78105986,0.02270522){\color[rgb]{0,0,0}\makebox(0,0)[lb]{\smash{$F$}}}%
  \end{picture}%
\endgroup%

%% file: allbeta.pdf_tex
\begingroup%
  \makeatletter%
  \providecommand\color[2][]{%
    \errmessage{(Inkscape) Color is used for the text in Inkscape, but the package 'color.sty' is not loaded}%
    \renewcommand\color[2][]{}%
  }%
  \providecommand\transparent[1]{%
    \errmessage{(Inkscape) Transparency is used (non-zero) for the text in Inkscape, but the package 'transparent.sty' is not loaded}%
    \renewcommand\transparent[1]{}%
  }%
  \providecommand\rotatebox[2]{#2}%
  \ifx\svgwidth\undefined%
    \setlength{\unitlength}{450.07998047bp}%
    \ifx\svgscale\undefined%
      \relax%
    \else%
      \setlength{\unitlength}{\unitlength * \real{\svgscale}}%
    \fi%
  \else%
    \setlength{\unitlength}{\svgwidth}%
  \fi%
  \global\let\svgwidth\undefined%
  \global\let\svgscale\undefined%
  \makeatother%
  \begin{picture}(1,0.21418415)%
    \put(0,0){\includegraphics[width=\unitlength]{allbeta.pdf}}%
    \put(0.938604,0.18133746){\color[rgb]{0,0,0}\makebox(0,0)[lb]{\smash{$\widetilde{R}^\alpha$}}}%
    \put(0.78170861,0.01297442){\color[rgb]{0,0,0}\makebox(0,0)[lb]{\smash{$F$}}}%
  \end{picture}%
\endgroup%

%% file: negligible.pdf_tex
\begingroup%
  \makeatletter%
  \providecommand\color[2][]{%
    \errmessage{(Inkscape) Color is used for the text in Inkscape, but the package 'color.sty' is not loaded}%
    \renewcommand\color[2][]{}%
  }%
  \providecommand\transparent[1]{%
    \errmessage{(Inkscape) Transparency is used (non-zero) for the text in Inkscape, but the package 'transparent.sty' is not loaded}%
    \renewcommand\transparent[1]{}%
  }%
  \providecommand\rotatebox[2]{#2}%
  \ifx\svgwidth\undefined%
    \setlength{\unitlength}{632.64003906bp}%
    \ifx\svgscale\undefined%
      \relax%
    \else%
      \setlength{\unitlength}{\unitlength * \real{\svgscale}}%
    \fi%
  \else%
    \setlength{\unitlength}{\svgwidth}%
  \fi%
  \global\let\svgwidth\undefined%
  \global\let\svgscale\undefined%
  \makeatother%
  \begin{picture}(1,0.4562468)%
    \put(0,0){\includegraphics[width=\unitlength]{negligible.pdf}}%
    \put(0.13724228,0.11101741){\color[rgb]{0,0,0}\makebox(0,0)[lb]{\smash{$E_\ep^k$}}}%
    \put(0.72738829,0.11121025){\color[rgb]{0,0,0}\makebox(0,0)[lb]{\smash{$E_\ep^{k+1}$}}}%
    \put(0.49125073,0.00770238){\color[rgb]{0,0,0}\makebox(0,0)[lb]{\smash{     \small{asymptotically negligible sets}}}}%
  \end{picture}%
\endgroup%

%% file: pinning2.pdf_tex
\begingroup%
  \makeatletter%
  \providecommand\color[2][]{%
    \errmessage{(Inkscape) Color is used for the text in Inkscape, but the package 'color.sty' is not loaded}%
    \renewcommand\color[2][]{}%
  }%
  \providecommand\transparent[1]{%
    \errmessage{(Inkscape) Transparency is used (non-zero) for the text in Inkscape, but the package 'transparent.sty' is not loaded}%
    \renewcommand\transparent[1]{}%
  }%
  \providecommand\rotatebox[2]{#2}%
  \ifx\svgwidth\undefined%
    \setlength{\unitlength}{320.4234375bp}%
    \ifx\svgscale\undefined%
      \relax%
    \else%
      \setlength{\unitlength}{\unitlength * \real{\svgscale}}%
    \fi%
  \else%
    \setlength{\unitlength}{\svgwidth}%
  \fi%
  \global\let\svgwidth\undefined%
  \global\let\svgscale\undefined%
  \makeatother%
  \begin{picture}(1,0.41445158)%
    \put(0,0){\includegraphics[width=\unitlength]{pinning2.pdf}}%
    \put(0.00142205,0.34626363){\color[rgb]{0,0,0}\makebox(0,0)[lb]{\smash{$A$}}}%
    \put(-0.00182365,0.20894433){\color[rgb]{0,0,0}\makebox(0,0)[lb]{\smash{$B$}}}%
  \end{picture}%
\endgroup%

%% file: pinning3.pdf_tex
\begingroup%
  \makeatletter%
  \providecommand\color[2][]{%
    \errmessage{(Inkscape) Color is used for the text in Inkscape, but the package 'color.sty' is not loaded}%
    \renewcommand\color[2][]{}%
  }%
  \providecommand\transparent[1]{%
    \errmessage{(Inkscape) Transparency is used (non-zero) for the text in Inkscape, but the package 'transparent.sty' is not loaded}%
    \renewcommand\transparent[1]{}%
  }%
  \providecommand\rotatebox[2]{#2}%
  \ifx\svgwidth\undefined%
    \setlength{\unitlength}{320.4234375bp}%
    \ifx\svgscale\undefined%
      \relax%
    \else%
      \setlength{\unitlength}{\unitlength * \real{\svgscale}}%
    \fi%
  \else%
    \setlength{\unitlength}{\svgwidth}%
  \fi%
  \global\let\svgwidth\undefined%
  \global\let\svgscale\undefined%
  \makeatother%
  \begin{picture}(1,0.41445158)%
    \put(0,0){\includegraphics[width=\unitlength]{pinning3.pdf}}%
    \put(-0.00229553,0.34020618){\color[rgb]{0,0,0}\makebox(0,0)[lb]{\smash{$A$}}}%
    \put(-0.00229553,0.06681739){\color[rgb]{0,0,0}\makebox(0,0)[lb]{\smash{$B$}}}%
  \end{picture}%
\endgroup%

%% file: cell2.pdf_tex
\begingroup%
  \makeatletter%
  \providecommand\color[2][]{%
    \errmessage{(Inkscape) Color is used for the text in Inkscape, but the package 'color.sty' is not loaded}%
    \renewcommand\color[2][]{}%
  }%
  \providecommand\transparent[1]{%
    \errmessage{(Inkscape) Transparency is used (non-zero) for the text in Inkscape, but the package 'transparent.sty' is not loaded}%
    \renewcommand\transparent[1]{}%
  }%
  \providecommand\rotatebox[2]{#2}%
  \ifx\svgwidth\undefined%
    \setlength{\unitlength}{189.93364258bp}%
    \ifx\svgscale\undefined%
      \relax%
    \else%
      \setlength{\unitlength}{\unitlength * \real{\svgscale}}%
    \fi%
  \else%
    \setlength{\unitlength}{\svgwidth}%
  \fi%
  \global\let\svgwidth\undefined%
  \global\let\svgscale\undefined%
  \makeatother%
  \begin{picture}(1,0.9505866)%
    \put(0,0){\includegraphics[width=\unitlength]{cell2.pdf}}%
    \put(0.46723934,0.00753346){\color[rgb]{0,0,0}\makebox(0,0)[lb]{\smash{$\delta_1$}}}%
    \put(-0.00194352,0.78147114){\color[rgb]{0,0,0}\makebox(0,0)[lb]{\smash{$\delta_2$}}}%
    \put(1.01901498,0.50844382){\color[rgb]{0,0,0}\makebox(0,0)[lb]{\smash{$2K\ep$}}}%
  \end{picture}%
\endgroup%